\documentclass[12pt]{amsart}
\usepackage{amssymb, amsmath, amsfonts}
\usepackage[raiselinks=false,colorlinks=true,citecolor=blue,urlcolor=blue,linkcolor=blue,bookmarksopen=true,pdftex]{hyperref}
\newtheorem{theorem}{Theorem}[section]
\newtheorem{proposition}[theorem]{Proposition}
\newtheorem{lemma}[theorem]{Lemma}
\newtheorem{remark}[theorem]{Remark}
\newtheorem{definition}[theorem]{Definition}
\newtheorem{example}[theorem]{Example}

\newcommand{\flag}{\textrm{Flag}}

\newcommand{\bz}{\mathbb{Z}}
\newcommand{\bq}{\mathbb{Q}}
\newcommand{\br}{\mathbb{R}}
\newcommand{\bc}{\mathbb{C}}

\newcommand{\bs}{\mathbb{S}}

\renewcommand{\hom}{\textrm{Hom}}
\newcommand{\wt}{\widetilde}

\newcommand{\pl}{\textrm{PL}}

\begin{document}
\baselineskip=15.5pt
\title[Formality of certain CW complexes]{Formality of certain CW complexes}  
\author{Prateep Chakraborty}
\author{Parameswaran Sankaran}
\address{The Institute of Mathematical Sciences, CIT
Campus, Taramani, Chennai 600113, India}
\email{prateepc@imsc.res.in}
\email{sankaran@imsc.res.in}


\date{}

\begin{abstract} 
Let $X$ be a simply connected path connected topological space which is formal in the sense of 
rational homotopy theory.  Let $Y=X\cup_\alpha\mathbb{D}^{n}$ where $\alpha:\mathbb{S}^{n-1}\to X$ 
is a non-torsion element.  Then we obtain a condition on $\alpha$ for the formality of $Y$.  We 
give several illustrative examples concerning the formality of a finite CW complex having only even dimensional 
cells.   

This is the corrected version of the earlier version which contained a serious error in Theorem 1.4. 
This theorem, which now Theorem 1.1 of this version, has now been corrected.  The proofs of Theorems 1.1, 1.2, and 1.3  
of the first version are not valid as they used the erroneous result.  In fact, we provide here a counterexample to 
the assertion of Theorem 1.1. (See Example 3.1 below.)   We do not know if the statement of Theorem 1.2, which asserted the formality of Schubert varieties in a generalized flag variety $G/B$, is valid.   Theorem 1.3 is correct as stated
as it had been proved previously by Panov and Ray using entirely different techniques.    
\end{abstract}
\maketitle
\section{Introduction}

Any path connected topological space $X$ 
has a functorial differential graded commutative algebra  (dgca) $A_\pl(X)$ over $\mathbb{Q}$, a {\it minimal model}  $(\mathcal{M}_X,d)$ (which 
is a dgca) and a 
dgc algebra morphism $\rho_X:\mathcal{M}_X\to A_\pl(X)$ such that $ \rho_X$ induces isomorphism in cohomology.  The minimal model is unique up to isomorphism.  The space $X$ is called {\it formal} if there exists a dgca morphism $(\mathcal{M}_X,d)\to (H^*(X;\bq),0)$ which induces isomorphism in cohomology.  If $X$ is a simply-connected space, its rational homotopy type is determined by $\mathcal{M}_X$.  In case $X$ is formal, $\mathcal{M}_X$ is determined by $H^*(X;\mathbb{Q})$ and so the rational homotopy 
type of $X$ is a `formal consequence of its cohomology algebra.'  

Let $X$ be a simply-connected formal CW complex and let $\alpha:\mathbb{S}^{n-1}\to X$ represent an element 
in the kernel of the Hurewicz homomorphism $\eta:\pi_{n-1}(X)\to H_{n-1}(X;\mathbb{Q})$.  Let $Y=X\cup_\alpha e^n$. 
We have the inclusion map $j:Y\hookrightarrow (Y,X)$ and the characteristic map $(D^n,\mathbb{S}^{n-1})\to (Y,X)$. 
Then $j^*:H^n(Y,X;\mathbb{Z})\to H^n(Y;\mathbb{Z})$ maps the positive generator $\wt{u}\in H^{n}(Y,X;\mathbb{Z})\cong H^n(D,\mathbb{S}^{n-1})\cong\mathbb{Z}$ to a non-zero element $u$ in $H^n(Y;\mathbb{Z}).$
If $\alpha$ represents a torsion element in $\pi_{n-1}(X)$, then $u$ is a non-zero indecomposable element in $H^n(Y;\mathbb{Q})$.  In this case, denoting the rationalization of $X$ by $X_0$, we see that $X_0\cup_\alpha e^n$ is homotopically equivalent to $X_0\vee \mathbb{S}^n$.  It follows that $Y$ is rationally equivalent to $X\vee \mathbb{S}^n$ which is a formal space.   

Our main result is the following.    Recall that a minimal model of a simply-connected space is isomorphic as a graded algebra to $\Lambda V$ 
where $V$ is a graded $\bq$-vector space $V=\oplus_{k\ge 2}V^k$ and $\Lambda V$ stands for the free graded-commutative algebra over $V$. Thus, $\Lambda V$ is isomorphic to the tensor product over $\bq$ of the 
symmetric algebra of $V^{even}=\oplus V^{2k}$ and the exterior algebra of $V^{odd}=\oplus V^{2k-1}$.  
One has the notion of lower gradation on $V$.  We recall these in \S2 as well the notion of standardness of  
lower gradation. 

\begin{theorem} \label{cell} 
Suppose that $X$ is a simply connected CW complex and is formal. Let $\mathcal{M}_X=\Lambda(V)$ and  suppose that $V=\oplus_{k\ge 0}V_k$ is a standard lower gradation of $V$.  Let 
$Y=X\cup_\alpha e^n$.   Suppose that $\eta([\alpha])=0$ so that $j^*(\wt{u})=:u\ne 0$.  
(i) If $[\alpha]\in \pi_{n-1}(X)$ is a torsion element then $u$ is indecomposable and $Y$ is formal.  
(ii) Let $[\alpha]\neq 0$ in $\pi_{n-1}(X)\otimes \mathbb{Q}$.   Suppose that 
$\langle v,[\alpha]\rangle =0$ for all $v\in V_k, k\ne 1,$  and that  $u$ is 
decomposable in $H^*(Y;\bq)$.   Then $Y$ is formal.  (iii) If $[\alpha]\in \pi_{n-1}(X)$ is not a torsion element and 
$u$ is not decomposable, then $Y$ is not formal.
\end{theorem}

{\it Throughout this paper $H^*(X)$ denotes the singular cohomology of $X$ with $\bq$-coefficients.  All differential graded commutative algebras will be over 
$\bq$.}


\section{Minimal models and formality}

In this section we recall the notion of a Sullivan algebra, 
a model for a cell attachment, and the 
stepwise construction of the minimal Sullivan model for 
a differential graded commutative cochain algebra with zero differential. We also 
prove Theorem \ref{cell}.  The reader is referred to \cite{fht} for a comprehensive 
treatment of rational homotopy theory.

\subsection{Sullivan algebra} \label{sullivanalgebra}
A differential graded commutative algebra (abbreviated {\it dgca}) $(M,d)$ 
is called a {\it Sullivan algebra} if the  following hold:
(i) {\it Freeness}: There exists a graded $\bq$-vector space $V=\oplus_{q\geq 1} V^q$ such that $M$ is freely generated by $V$, that is, 
$M=\Lambda V:=S^*(V^{even})\otimes E^*(V^{odd})$ 
 where $V^{even}=\oplus_{q\geq 1} V^{2q}, ~V^{odd}=\oplus_{q\geq 1} V^{2q-1}$. Here $S^*(V)$ denotes the symmetric algebra on $V$ 
and $E^*(V)$ denotes the exterior algebra on $V$.  
(ii) {\it  Nilpotence:} There is a well-ordering on a basis 
$\{v_\alpha\}$ 
of $V$ consisting of homogeneous elements such that 
for each $\alpha$, $d(v_\alpha)$ is a polynomial in the $v_\beta, \beta<\alpha.$   

The nilpotence condition can be restated as follows:
There is an increasing filtration $V=\cup_{k\geq 0} V(k)$, such that $d(V(k))\subset \Lambda(V(k))$ and there 
exists a subspace $V_k\subset V(k)$ such that $d(V_k)\subset \Lambda(V(k-1))$  and $\Lambda(V(k))=\Lambda(V_k)\otimes \Lambda(V(k-1))$.  This filtration is 
referred to as {\it the lower filtration} of $M$.

A {\it Sullivan model} $(M,d)$ for a dgca $(A,d)$ is a Sullivan 
algebra $(M,d)$ together with a dgca morphism 
$f: M\to A$ which is a quasi-isomorphism. Thus, $f$ is a dgca morphism which induces 
an isomorphism $f^*:H^*(M,d)\to H^*(A,d)$.  A Sullivan 
model $(M,d)$ is called a {\it minimal model} if $d(M)\subset M^+.M^+$, the ideal of decomposable elements.  A dgc algebra $A$ with $H^0(A)=\bq$ has a unique minimal model up to isomorphism. Such a dgc algebra is called {\it formal} if there exists a dgca morphism $\Phi:(\mathcal{M}_A,d)\to (H^*(A),0)$ which is a quasi-isomorphism where $\mathcal{M}_A$ denotes the minimal model of $A$. 

Suppose that $X$ is a path connected topological space.  Sullivan \cite{sullivan} constructed a natural dgc algebra 
$(A_\pl(X),d)$ over $\bq$ called {\it the polynomial differential forms on $X$} over $\bq$, which is contravariant in $X$. Its cohomology  $H^*(A_\pl(X),d)$ is naturally isomorphic to 
$H^*(X;\bq)$. See \cite[\S 10]{fht} for 
details of construction of the functor $A_\pl(-)$. 
   
A Sullivan model (resp. minimal model) for $X$ is by definition a Sullivan model (resp. minimal model) for $A_\pl(X)$.  Any path connected space $X$ has a 
minimal Sullivan model $\mathcal{M}_X$.  See \cite[Proposition 12.1]{fht}.  Minimal models of $X$ are unique up to isomorphism provided $H^1(X)=0$.  If $X$ is simply connected, then $\hom(\pi_k(X), \bq)\cong V^k$ where $\mathcal{M}_X=\Lambda(V), V=\oplus_{k\ge 2} V^k$ denotes the minimal model of $X$. If $X$ and $Y$ are simply connected and have the same rational homotopy type, then their minimal models are isomorphic (as dgc algebras).  In fact we have 
a bijection between the collection of all the rational homotopy types of 
simply connected spaces and the collection of  all isomorphism classes of 
minimal Sullivan algebras over $\bq$.   Assume that $X$ and $Y$ are simply connected and that their 
rational cohomology algebras are of finite type. Then we have an isomorphism of sets: $[X_0,Y_0]\to [\mathcal{M}_Y,\mathcal{M}_X]$ where $X_0$ denotes the rationalization of $X$, $[A,B]$ denotes the 
homotopy classes of dgca morphisms $A\to B$ between Sullivan algebras
and $\mathcal{M}_X:=\mathcal{M}_{A_\pl(X)}$ is the minimal model of $X$. Observe that $\mathcal{M}_X$ and $\mathcal{M}_{X_0}$ are naturally isomorphic 
since $X\subset X_0$ is a rational homotopy equivalence. 
The isomorphism is obtained by 
sending $[f]\in[X_0,Y_0]$ to the homotopy class of 
any lift $\theta:\mathcal{M}_Y\to \mathcal{M}_X$ 
of $A_\pl(f):A_\pl(Y_0)\to A_\pl(X_0)$ so that the following diagram commutes up to homotopy:
\[\begin{array}{rll} 
\mathcal{M}_Y &\stackrel{\theta}
{\longrightarrow}&\mathcal{M}_X\\
\rho_Y\downarrow& & \downarrow \rho_X\\
A_\pl(Y_0)&\stackrel{A_\pl(f)}{\longrightarrow}&A_\pl(X_0).\\
\end{array}
\]
\noindent
{\it Notations.} 
If $V$ is a graded vector space, then $V^{\leq k}$ (resp. $V^{<k}$) denotes the subspace consisting of elements of degree atmost $k$ (resp. $k$). If $A$ is a differential graded algebra, $A^{\le q}$ (resp. $A^{<k}$) denotes the differential graded subalgebra of $A$ generated by elements of degree atmost (resp. less than) $k$.

If $A$ is a dgc algebra, we denote by $\mathcal{D}(A)$ (or simply $\mathcal{D}$ if $A$ is clear from the context) the ideal of decomposable elements in $A$. By abuse of notation we write $A^n/\mathcal{D}$ to mean $A^n/\mathcal{D}\cap A^n\subset A/\mathcal{D}$.

Two dgc algebras $(A,d)$ and $(B,d)$ are quasi-isomorphic if there is a finite sequence of dgca  morphisms $f:=\{f_i\}$ where $A_0\stackrel{f_0}{\to}A_1\stackrel{f_1}{\leftarrow} A_2\stackrel{f_2}{\to} \cdots\stackrel{f_{2n-1}}{ \leftarrow}A_{2n}$ with $(A_0,d)=(A, d), (A_{2n},d)=(B,d)$ such that induced 
morphisms in cohomology are all isomorphisms. In this 
case we write $(A,d)\stackrel{f}{ \leftrightarrow}(B,d)$ or $(A,d)\simeq (B,d)$. We denote by $f^*:H^*(A,d)\to H^*(B,d)$ the composition of isomorphisms $(f_{2n-1}^*)^{-1}\cdots \circ f_0^*$.

\subsection{Minimal model of $(A,0)$} \label{minimalmodel}
We refer the reader to \cite{fht} for  construction of the minimal model $(\mathcal{M}_A,d)$ 
for a dgca $(A,d)$.   For our purposes, we need only consider minimal model for a dgca with zero differential satisfying 
$A^0=H^0(A)=\mathbb{Q}, A^1=H^1(A)=0$.   We shall 
particularly use the description given in \cite[\S3]{hs}. 

\begin{lemma}\label{nolambda0}
Let $(\mathcal{M}_A,d)=(\Lambda V,d)$ be a minimal model of a dgca $(A,0)$ with zero differential.  Let $\rho_A$ (or more briefly $\rho$) denote a quasi-isomorphism $\mathcal{M}_A\to A$ inducing identity in cohomology.  Then 
there exists a lower gradation $V=\oplus_{k\ge 0}V_k$ such that (i) $\rho(V_k)=0$ for all $k\ge 1$, and, 
(ii) $d(V_k)\subset \Lambda(V_0).\Lambda^+(V_1\oplus\cdots\oplus V_{k-1})$ for $k\ge 2$. 
\end{lemma} 
\begin{proof} The existence of a lower gradation $V_k,k\ge 1,$ such that $\rho(V_k)=0$ is well-known.  Indeed (i) holds by 
the construction of the minimal model of $A$ given in \cite[\S3]{hs}.  We start with such a lower gradation $V_k, k\ge 0$ and modify this to obtain a new lower gradation $V'_k$ so as to meet both our requirements.  We set $V_k^n=V_k\cap V^n$.

Let $\{y_\gamma\}_{\gamma\in J_{k,2}}$ be a basis for $V_2^k$.  Write $dy_\gamma=u_0+u_1$ where $u_0\in \Lambda(V_0)^{k+1}$ and $u_1\in \Lambda(V_0).\Lambda^+(V_1)$.  Then $\rho(u_1)=0$ using $\rho(V_1)=0$ and the fact that $\rho$ is an algebra homomorphism.  Therefore, $0=d\rho(y_\gamma) 
=\rho(dy_\gamma)=\rho(u_0)$ implies that $u_0=\sum f_i.dv_i=d(\sum f_iv_i)$ where $f_i\in \Lambda(V_0), v_i\in V_1$ since $u_0\in \Lambda(V_0)$ and $\rho$ induces isomorphism in cohomology.  Now let $y_\gamma'=y_\gamma-\sum f_iv_i.$  Then $dy_\gamma'=u_1\in 
\Lambda(V_0).\Lambda^+(V_1)$ and $\rho(y_\gamma')=\rho(y_\gamma)-\sum \rho(f_i)\rho(v_i)=0$ as $\rho(V_2)=0=\rho(V_1)$.   We define $V_2'^k\subset \Lambda(V_0\oplus V_1)^k\oplus V_2^k$ to be the space spanned by $y_\gamma', \gamma\in J_{k,2}$.  Set $V_2'=\oplus_{k\ge 3}V^{'k}_2$.  Note that $V_2'\cap(V_0\oplus V_1)=0$, $ V(2)=V_0+ V_1+V_2'$, $\rho(V_2')=0$ and $d(V_2')\subset \Lambda(V_0).\Lambda^+(V_1)$.   
 
We now proceed by induction.  Assume that $V_j', 2\le j<n,$ have been constructed satisfying (i) and (ii) such 
that $V_0+V_1+V_2'+\cdots+V_{n-1}'=V(n-1)$. It is convenient to set $V_1':=V_1.$
Let $\{y_\gamma\}_{\gamma\in J_{k,n}}$ be a basis for $V_n^k$.  Write $dy_\gamma=z_0+z_1$ where $z_0\in \Lambda(V_0)^k$ and $z_1\in \Lambda(V_0).\Lambda^+(V_1\oplus V_2'\oplus \cdots\oplus V'_{n-1})$.  Then 
$\rho(z_1)=0$ using $\rho(V'_j)=0, j\ge 1$.  Therefore, $0=d\rho(y_\gamma)
=\rho(dy_\gamma)=\rho(z_0)$ implies that $z_0=d(\sum f_j.x_j)$ where $f_j\in \Lambda(V_0), x_j\in \Lambda^+(V_1)$ since $z_0\in \Lambda(V_0)$.  Set $y'_\gamma:=y_\gamma-\sum f_jx_j$.  Then  $dy'_\gamma=z_1\in \Lambda(V_0).\Lambda^+(V_1\oplus \cdots\oplus V_{n-1}')$ and $\rho(y_\gamma')=0$ as $\rho(y_\gamma)=0$ and $\rho(V_j)=0, 1\le j<n$.    Thus $V'_n:=\oplus_{k\ge 3}(\oplus_{\gamma\in J_{k,n}} \mathbb{Q}y'_{\gamma}) $ satisfies (i) and (ii).  Furthermore $V(n)=V(n-1)+V'_n$, $V_n\cap V(n-1)=0$.
This completes the induction step and we  
see that $V_0,V_1, V_j', j\ge 2,$ yield a lower gradation for $V$ that meets our requirements. 
  \end{proof}

\begin{definition}\label{standard} Let $\mathcal{M}_A$ be a minimal model of $(A,0)$. 
We say that a lower gradation $V=\oplus_{k\ge 0} V_k$ of $\mathcal{M}_A=\Lambda(V)$ is {\em standard} if 
it satisfies conditions (i) and (ii) of Lemma \ref{nolambda0}.
\end{definition}

\subsection{A model for cell attachment}\label{cellmodel} Let $X$ be a simply connected topological space. 
Let $Y=X\cup_\alpha e^n$ where $\alpha:\bs^{n-1}\to X$ represents an element $[\alpha]\in \pi_{n-1}(X)$. We assume that $n\geq 2$ so that $Y$ is also simply connected.
We recall the following proposition which will 
play a crucial role in our proofs.   Let $m_X:(\mathcal{M}_X,d)\to (A_\pl(X),d)$ be a minimal Sullivan model 
for $X$.  Suppose that $\mathcal{M}_X=\Lambda(V)$ so 
that $V=\oplus_{ k\geq 2} V^k$.  (Note that $V^1=0$ since $X$ is simply connected.)  Recall that $V^k\cong \hom(\pi_k(X),\bq)$;
thus we have the pairing $\langle -, -\rangle:V^k\times \pi_k^\bq(X)\to \bq$ defined by evaluation.   If $n=2$, then $[\alpha]=0$ 
as $X$ is simply connected.  It follows that $Y\simeq X\vee \bs^2$ which is formal.  

Let $n\geq 3$. 
Let $M_\alpha=\Lambda(V_\alpha)$ be the dgca defined as follows: $V_\alpha:=V\oplus \bq u_\alpha$,  $\deg(u_\alpha)=n$, $u_\alpha^2=u_\alpha.v=0, v\in V$, with differential $d_\alpha$ where $d_\alpha(u_\alpha)=0$ and $d_\alpha(v)=
d v+\langle v, \alpha\rangle u_\alpha,\  \ v\in V$.  

\begin{proposition}  \label{model}
The dgca $(M_\alpha, d_\alpha)$ defined above is a model for 
$Y=X\cup_\alpha e^n$.  Moreover, 
one has the following diagram of dgc algebras in which the rows are exact and the vertical arrows are quasi-isomorphism:
\[
\begin{array}{ccccccc}
0\to &\bq u_\alpha &\hookrightarrow &M_\alpha & \stackrel{\lambda}{\to}& \mathcal{M}_X&\to 0\\
&\updownarrow\simeq & & \updownarrow\mu & & \downarrow \Phi &\\
0 \to & A_\pl(Y,X)& \stackrel{A_\pl(j)}{ \to} &A_\pl(Y) &\stackrel{A_\pl(i)}{\to} & A_\pl(X)  &\to 0\\
&\downarrow\simeq &&&&&\\
& A_\pl(D^n, \bs^{n-1}) & &&&&
\end{array}
\]
where $i:X\hookrightarrow Y$ and $j:Y\hookrightarrow (Y,X)$ are inclusions and $\lambda$ is induced by projection $V_\alpha\to V$.  The induced diagram 

\[
\begin{array}{ccccccc}
0\to &\bq u_\alpha &\hookrightarrow&H^*(M_\alpha) & \stackrel{\lambda^*}{\to}& H^*(\mathcal{M}_X)&\to 0\\
&\downarrow\cong & & \downarrow\mu^* & & \downarrow \Phi ^*&\\
0 \to & H^*(Y,X)& \stackrel{j^*}{ \to} &H^*(Y) &\stackrel{i^*}{\to} & H^*(X)  &\to 0\\
&\downarrow\cong &&&&&\\
& H^*(D^n, \bs^{n-1}) & &&&&
\end{array}\eqno(1)
\]
is commutative with exact rows in which the vertical arrows are all isomorphisms.  \hfill $\Box$
\end{proposition}
We refer the reader to \cite[Chapter 13]{fht} for a proof.

\begin{remark}
The dgca $M_\alpha$ is not a minimal model 
for $Y$ most often. Indeed it is not free except in the case $V=0$ and $n$ odd, since $u_\alpha^2=0$ and the relation $u_\alpha.v=0$ holds for $v\in V$.      
\end{remark}   
    

\subsection{Proof of Theorem \ref{cell}}
We keep the notations and set-up of \S \ref{cellmodel}.  It is understood that 
a base point for $X$ is chosen and fixed and it serves as the point for $Y=X\cup_\alpha e^n$ as well; the homotopy 
groups are defined with respect to this choice and 
will be suppressed in the notation $\pi_k(X)$, etc.
Recall that $i$ (resp. $j$) denotes the inclusion $X\hookrightarrow Y$ (resp. $Y\hookrightarrow (Y,X)$).  Also $V$ and $W$ are graded vector spaces so that $\mathcal{M}_X=\Lambda(V)$ 
and $\mathcal{M}_Y=\Lambda(W)$.  One has a morphism 
of dgca $\phi:\mathcal{M}_Y\to \mathcal{M}_X$ which is a 
lift of $A_\pl(i): A_\pl(Y)\to A_\pl(X)$.  The 
linear part $Q(\phi):W\to V$ of $\phi$ is defined by 
the requirement that $\phi(w)-Q(\phi(w))\in \Lambda^{\geq 2}V$; it induces $i^*: \hom(\pi_k(Y),\bq)\to \hom(\pi_k(X),\bq)$ for all $k$ 
under the isomorphisms $V^k\cong \hom(\pi_k(X),\bq)$ 
and $W^k\cong \hom(\pi_k(Y),\bq)$.

Recall that $X$ is simply connected.
By the relative Hurewicz theorem, we obtain that $\eta: \pi_n(Y,X)\cong H_n(Y,X;\bz)\cong \bz$. 
The group $\pi_n(Y,X)=\bz$ is 
generated by the homotopy class of the characteristic map $ \wt{\alpha}: (D^n,\bs^{n-1})\to (Y,X)$ of the cell $e_\alpha$.  
The homomorphism $\partial: \pi_n(Y,X)\to \pi_{n-1}(X)$ maps $[\wt{\alpha}]$ to $[\alpha]$. Denoting by $\pi_k^\bq$ the rational homotopy group functor $\pi_k(-)\otimes \bq$, we have the commuting diagram 

\[\begin{array}{lllllll}
&&\pi_n^\bq(Y)&\stackrel{j_*}{\to}&\pi_n^\bq(Y,X)&\stackrel{\partial}{\to} &\pi_{n-1}^\bq(X) \\ 
&&\eta\downarrow &&\eta\downarrow &&\downarrow\eta\\
H_{n}(X)&\to &H_{n}(Y)& \stackrel{j_*}{\to} & H_{n}(Y,X)&\stackrel{\partial}{\to}& H_{n-1}(X).\\ 
\end{array}\eqno(2)
\]
where $\eta$ denotes the Hurewicz homomorphism,    
with the middle one being an isomorphism. 

Suppose that $\eta([\alpha])\neq 0$.  Then $\partial (\eta[\wt{\alpha}])=\eta([\alpha])\neq 0$ and since $H_n(Y,X)\cong 
\bq$, we conclude that $j_*=0$ and $i_*:H_n(X)\to H_n(Y)$ is an 
isomorphism since $H_{n+1}(Y,X)=0$.  Therefore 
$i^*:H^{n}(Y)\to H^n(X)$ is also an isomorphism.

Suppose that $\eta([\alpha])=0$. 
This happens, for example, when $\alpha$ is a torsion element or when $H_{n-1}(X)=0$. 
(However $\eta[\alpha]=0$ does not imply that $\alpha$ is of finite order. For example, one can choose $\alpha$ to be an element of infinite order in $\pi_{4m-1}(\bs^{2m})$.)
There exists an element $\gamma\in H_n(Y)$ such that $j_*(\gamma)
=\eta([\wt{\alpha}])$. Let $\wt{u}$ denote the generator of $H^n(Y,X)=\hom(H_n(Y,X),\bq)\cong\bq$ such that $\langle \wt{u}, \eta([\wt{\alpha}])\rangle =1$. Then $j^*(\wt{u})=:u$ 
is a non-zero element of $H^n(Y)$ and we have 
$\langle u, \gamma\rangle =1$. 
Therefore, using the 
exact sequence $H^{n-1}(Y,X)\stackrel{j^*}{\to} H^{n-1}(Y)\stackrel{i^*}{\to} H^{n-1}(X)\to H^n(Y,X)\stackrel{j^*}{\to} H^n(Y)\stackrel{i^*}{\to} H^n(X)\to H^{n+1}(Y,X)=0$, we have 
\[H^n(Y)\cong\left\{ \begin{array}{ll}
H^n(X)\oplus \bq u &\textrm{if} \  \ \eta([\alpha])=0,\\
H^n(X) & \textrm{if} \  \ \eta([\alpha])\neq 0,
\end{array}
\right . \eqno(3) 
\]
and 
\[ H^{n-1}(X)\cong \left\{ \begin{array}{ll}
H^{n-1}(Y) & \textrm{if} \ \ \eta([\alpha])=0,\\
H^{n-1}(Y)\oplus \bq \wt{u}& \textrm{if} \ \ \eta([\alpha])\ne 0.
\end{array}\right . \eqno(4)\] 
Since $\hom(\pi_{n-1}(X),\bq)\cong V^{n-1}$, using the exactness of the sequence $\pi_n^\bq(Y,X)\stackrel{\partial}{\to} \pi_{n-1}^\bq(X)\stackrel{i_*}{\to} \pi_{n-1}^\bq(Y)\to \pi_{n-1}^\bq(Y,X)=0$ we 
see that 
\[V^{n-1}\cong\left\{\begin{array}{ll}
W^{n-1}\oplus\bq & \textrm{if}\  \ [\alpha]\neq 0,\\
W^{n-1}& \textrm{if}\  \  [\alpha]= 0,
\end{array}\right .  \eqno(5)
\]
via the restriction of $Q(\phi)$. 

Summarizing the above discussion we obtain the 
following.

\begin{lemma} \label{cohomologyofy}
(i) Suppose that $\eta[\alpha]=0$.  Then $j^*(\wt{u})=u\ne 0$ 
and $H^n(Y)\cong H^n(X)\oplus \bq u$,  $H^{k}(Y)\cong H^{k}(X), k\ne n$.  If $[\alpha]\ne 0$ in $\pi^\bq_{n-1}(X)$, then $V^{n-1}\cong W^{n-1}\oplus \bq.$\\
(ii) Suppose that $\eta[\alpha]\ne 0$.  Then $j^*(\wt{u})
=0$ and $H^k(Y)\cong H^k(X), k\ne n-1$,  $H^{n-1}(X)\cong H^{n-1}(Y)\oplus \bq [\wt{u}]$. Moreover $V^{n-1}\cong W^{n-1}\oplus \bq$.\hfill $\Box$
\end{lemma}

We now establish Theorem \ref{cell}.

\noindent
\textit{Proof of Theorem \ref{cell}.} 
(i) If $[\alpha]\in \pi_{n-1}(X)$ is a torsion element then $Y$ is rational homotopically equivalent to $X_0\vee \bs^{n}$. Hence $Y$ is formal. 

(ii) In this case $Q(\phi):W^k\to V^k$ is an isomorphism for $k\le n-2$ and is a monomorphism when $k=n-1$.  Moreover, $Q(\phi)(W^{n-1})=\ker([\alpha])\subset V^{n-1}$ has codimension $1$. 
Write
$Q(\phi)(W^{n-1})\oplus \mathbb{Q}v_\alpha=V^{n-1}$ where $v_\alpha\in V^{n-1}$ 
is an element such that $\langle v_{\alpha},[\alpha]\rangle=1$.  (We shall presently make a more specific 
choice of $v_\alpha$.)  
Write 
$u=P(\overline{v}_1,\ldots, \overline{v}_r)$ with $\overline{v}_q\in H^{<(n-1)}(Y)\cong 
H^{<(n-1)}(X)$ where $\overline{v}_q$ are indecomposable elements. 
Since $\Lambda(V_0)\to H^*(X)$ is onto, we choose cocycles $v_q\in V_0\subset \mathcal{M}_X$ so that $v_q\mapsto \overline{v}_q$. 
We set $w=P(v_1,\ldots,v_r)\in \mathcal{M}_X$. 
Since $X$ is formal we have a quasi-isomorphism $\Phi:(\mathcal{M}_X,d)\to (H^*(X),0)$. Since $i^*(u)=0$ we have $\Phi(P(v_1,\ldots, v_r))=0$ in $H^n(X)$. That is, $P(v_1,\ldots,v_r)=:w\in \ker(\Phi)$.  Note that $w\in \mathcal{M}_X$ is a cocyle. 
Since $\Phi^*$ is a monomorphism, $w=d_X(v_\alpha)$ for some $v_\alpha\in \Lambda V(1)\subset \mathcal{M}_X^{n-1}$. 
We claim that $\langle v_\alpha, [\alpha]\rangle \neq 0$. 
Indeed, since $\mu:(M_\alpha,d_\alpha)\leftrightarrow A_{\pl}(Y)$ is a quasi-isomorphism we have, using the commutative diagram (1) of \cite{cs}, that $\mu^*([w])=P(\overline{v}_1,\ldots,\overline{v}_r)=u\in H^n(Y)$ is non-zero.  So $[w]\neq 0$ in $H^n(M_\alpha).$  If $\langle v_\alpha, [\alpha]\rangle=0$, then  
$d_\alpha(v_\alpha)=d_X(v_\alpha)=w$, whence $\mu^*([w])=0$ in $H^n(M_\alpha)$, a contradiction.  Therefore $\langle v_ \alpha,[\alpha]\rangle \neq 0$.  Now this implies 
that $v_\alpha\notin Q(\phi)(W)=\ker ([\alpha]).$

By hypothesis $\langle v,[\alpha]\rangle=0~ \forall v\in V^{n-1}_0\oplus (\oplus_{k\ge 2}V_k^{n-1})$.  

The surjective homomorphism $i^*:H^n(Y)\to H^n(X)$ induces an isomorphism $H^n(Y)/\mathcal{D}\to H^n(X)/\mathcal{D}\cong V_0^n$.  (Here $\mathcal{D}$ stands for the space of decomposable elements.) Choose a linear map $\theta':V_0^n\to H^n(Y)$ such that 
$\Phi(v)=i^*(\theta'(v)), v\in V_0^n$ and extend it to a linear map $\theta:V^n\to H^n(Y)$ by setting 
$\theta(v)=0$ for $v\in \oplus_{k\ge 1}V_k^n$.    

Define a vector space homomorphism $\psi:V\oplus \mathbb{Q}u_\alpha \to H^*(Y)$ as follows:
\[\psi(v)=\left \{ \begin{array}{ll} (i^*)^{-1}\circ \Phi(v)&  \textrm{if} \  \ v\in V^k, \ \ k\ne n,\\
\theta (\Phi(v)) & \textrm{if} \ \ v\in V^{n},\\ 
-u/\langle v_\alpha,[\alpha]\rangle & \textrm{if} \ \ v=u_\alpha.
\end{array}
\right .
\]
This extends to a homomorphism $M_\alpha\to H^*(Y)$, again denoted by  $\psi$, of the graded commutative algebra $M_\alpha$ because the relations $u^2=0, u.z=0$
for all $z\in H^*(Y)$ hold.  
Note that $\psi(w)=u$.

We claim that $\psi$ is a dgca morphism, that is, $\psi\circ d_\alpha=0$.   
Clearly $\psi(d_\alpha(u_\alpha))=\psi(0)=0$. 
Let $v\in V^k, k\ne n-1$.  Then $d_\alpha(v)=d_Xv$ and so $\psi(d_\alpha(v))=(i^*)^{-1}(\Phi(d_Xv))=0$.
If $v\in V^{n-1}_j, j\ne 1$, then $d_\alpha(v)=d_Xv$ since $v\in \ker([\alpha])$. If $j=0$, then $d_X(v)=0$ whence $\psi(d_\alpha v)=0$. Assume that $j>1$. Since $V=\oplus_{k\ge 0} V^k$ is a standard lower gradation, we see that $d_Xv$ is a sum of monomials in each of which {\it  there is a factor belonging to $V_i, 1\le i<j,$ present} by Lemma  \ref{nolambda0}.  Therefore $\psi(d_Xv)=0$.

Finally, let $v\in V_1^{n-1}=V_1^{n-1}\cap \ker([\alpha])\oplus \mathbb{Q}v_\alpha$. Suppose $v\in \ker([\alpha]).$  Then $d_\alpha(v)=d_Xv=f(w_1,\ldots, w_s)$.  
Since $\mu^*: H^*(M_\alpha)\to H^*(Y)$ is an isomorphism of algebras which agrees with $H^*(\mathcal{M}_X)\to H^*(X)$ in degrees less than 
$n-1$, we see that $0=[d_\alpha v]$ in $H^*(M_\alpha)$ 
implies that $[d_Xv]=0$ in $H^n(Y)$. On the 
other hand $\psi(d_Xv)=\psi(f(w_1,\ldots, w_s))
=f(\psi(w_1),\ldots, \psi(w_s))$ is the image of 
the element $f(w_1,\ldots, w_s)$ under $\mu^*$.  As 
$f(w_1,\ldots, w_s)=d_\alpha v$, we conclude that 
$\psi(d_\alpha v)=0$.  

It remains to consider the case $v=v_\alpha$. Then $d_\alpha(v_\alpha)=d_Xv_\alpha+\langle v_\alpha, [\alpha]\rangle u_\alpha=
w+\langle v_\alpha, [\alpha]\rangle u$.   It follows that $\psi(d_\alpha v_\alpha)=\psi(w)-u=0$.     
It is clear that $\psi$ induces isomorphism in cohomology.

Since $M_\alpha\simeq \mathcal{M}_Y,$ there exists a dgca morphism $h:\mathcal{M}_Y\to M_\alpha$ 
which induces isomorphism in cohomology.   Then $\psi\circ h:\mathcal{M}_Y\to H^*(Y)$ induces isomorphism in cohomology.


(iii)  Let $[\alpha]\ne 0$ in $\pi_{n-1}^\bq(X)$.
Assume that $j^*(\wt{u})=u$ is {\it not} decomposable, and that $Y$ is formal. We shall arrive at 
a contradiction.  Recall that by Proposition \ref{model}, 
$\mu:M_\alpha\leftrightarrow A_\pl(Y)$ is a quasi-isomorphism.
 
Let $\nu: \mathcal{M}_Y\to M_\alpha$ 
and $\Psi:\mathcal{M}_Y\to H^*(Y)$ be quasi-isomorphisms so that $\mu^*\circ\nu^*=\Psi^*$.     
Let $\lambda:M_\alpha\to \mathcal{M}_X$ be the dgca 
morphism considered in Proposition \ref{model}.  Then 
$\phi=\lambda\circ \nu:\mathcal{M}_Y\to \mathcal{M}_X$ is a lift of $A_\pl(i): A_\pl(Y)\to A_\pl(X)$ and induces $i^*:H^*(Y)\to H^*(X)$.  From (5), we have an  
isomorphism $V^{n-1}\cong W^{n-1}\oplus \bq$ given by $Q(\phi)$.  Also $V^k\cong W^k$ if $k\leq n-2$.  

Consider the dgca morphism $\iota: \Lambda(W^{\leq n-1})\to M_\alpha$.  Then $\iota^*$ is an 
isomorphism in dimension $\le n-2$ and the cokernel of $\iota^*:H^n(\Lambda(W^{\le n-1}))\to H^n(M_\alpha)\cong H^n(Y)$ is isomorphic to $\ker(d_Y)\cap W^n$. Since $[dv_\alpha]=u$ and since $dv_\alpha \in \Lambda(V^{\le n-2})=\Lambda(W^{\le n-2})$, 
we see that $u$ belongs to the image of $(\iota^*)$. 
Since $u$ is indecomposable we see that 
$\dim(coker(\iota^*))<\dim  (H^n(M_\alpha)/\mathcal{D})=\dim (H^n(Y)/\mathcal{D})$.
Therefore $\dim(\ker(d_Y)\cap W^n)<\dim (H^n(Y)/\mathcal{D})$. 
On the other hand, since $Y$ is formal, 
$\dim (H^n(Y)/\mathcal{D})=\dim(\ker(d_Y)\cap W^n)$.
Therefore we conclude that $Y$ cannot be formal. 
\hfill $\Box$

\section{Examples}

In this section we construct various illustrative examples.   

In view of the above Theorem \ref{cell}, we call $\alpha:\mathbb{S}^{n-1}\to X$  or the element $[\alpha]\in \pi_{n-1} (X)\otimes \mathbb{Q}$ {\it special} 
if $[\alpha]$ is in the kernel of the Hurewicz homomorphism and $\langle v,[\alpha]\rangle=0$ for all $v\in V_{k}\cap V^{n-1}, k\ne 1$.      

Our first example shows that merely assuming $u$ to be 
decomposable is not sufficient to conclude formality of $Y$ as claimed in Theorem 1.4 of \cite{cs}.  

\begin{example} {\em 
(1) Let $X=X_0\vee X_1$ where $X_0=\mathbb{S}^2\vee \mathbb{S}^2\vee \mathbb{S}^2, X_1=\{(x,y,z)\in \mathbb{S}^2\times \mathbb{S}^2\times \mathbb{S}^2\mid \textrm{~ $x$~or~$y$ ~or~ $z$ equals ~} *\}$ where $*$ denotes the base point 
of $\mathbb{S}^2$.  Then $X_0$ is formal, being a wedge of formal spaces.  The space $X_1$ is formal 
since it is the $4$-skeleton of the formal $\mathbb{S}^2\times \mathbb{S}^2\times \mathbb{S}^2$ where the $2$-sphere 
is given the cell structure with one $0$-cell and one $2$-cell.  Therefore $X$ is formal.  One has $H^*(X_0;\mathbb{Q})
=\mathbb{Q}[a_1,a_2,a_3]/\langle a_1^2,a_2^2,a_3^2,a_1a_2,a_2a_3,a_3a_1\rangle, |a_i|=2$ and $H^*(X_1;\mathbb{Q})
=\mathbb{Q}[x_1,x_2,x_3]/\langle x_1^2,x_2^2,x_3^2,x_1x_2x_3\rangle, |x_i|=2$.  The cohomology algebra of 
$X$ is readily computed from this.   The minimal model $(\mathcal{M}_X,d)=(\Lambda(V),d)$ of $X$  can be computed 
from the description of $H^*(X;\mathbb{Q})$ since $X$ is formal.  We obtain that $V_0=V^2=H^2(X;\mathbb{Q})$ is six dimensional, with basis 
$a_1,a_2,a_3,x_1,x_2,x_3$ all in degree $2$ where $d(a_i)=d(x_i)=0$,  $V^3\subset V_1$ has basis $c_{ij}, v_{i}, w_{ij}, 1\le i,j\le 3$
where  $dc_{ij}=a_ia_j, dv_{i}=x_i^2, dw_{ij}=a_ix_j$.  Also note that there are elements $f_{ij}\in V^4, i\ne j$ such that $df_{ij}=a_ic_{ij}-a_jc_{ii}$, elements $g_{ij}\in V^5,i\ne j,$ such that $dg_{ij}=a_jf_{ij}-a_if_{ji}+c_{ii}c_{jj}$ 
and there is an element $z\in V_0, |z|=5,$ with $dz=x_1x_2x_3$.   Let $\alpha:\mathbb{S}^5\to X$ be such that 
$\langle g_{12}, [\alpha]\rangle =l, \langle z,[\alpha]\rangle=m$ where $l,m$ are non-zero integers. Then
$y:=mg_{12}-lz$ vanishes on $[\alpha]$, $d_\alpha y=mdg_{1,2}-lx_1x_2x_3, d_\alpha z=x_1x_2x_3+mu_\alpha$.  So, 
in $H^*(M_\alpha)\cong H^*(Y)$ we obtain $[u_\alpha]=u=(-1/m)x_1x_2x_3$, a decomposable element.  {\it Note that 
$g_{12}\in V_3\cap V^5$ and  $\langle -, [\alpha]\rangle :V_3\to \mathbb{Q}$ is non-zero.}

We claim that $Y$ is not formal.  Suppose that $Y$ is formal then it is readily seen that, writing $\mathcal{M}_Y =\Lambda W$, $W^j=V^j, j\le 4,$ and $V^5\cong W^5\oplus \mathbb{Q} z$ where $W^5$ is identified with the kernel of $[\alpha]$.  
In particular $y\in W^5$.  Let $\Psi:\mathcal{M}_Y\to H^*(Y)$ be a dgca morphism 
that induces isomorphism in cohomology.    Then $\Psi|W^j=\Phi|W^j, j\le 3,$ where $\Phi:\mathcal{M}_X\to H^*(X)$ is a suitable quasi-isomorphism.  Since $\Psi(a_i)=\Phi(a_i)=a_i$ and since $a_i.H^4(Y;\mathbb{Q})=0$ we get 
$0=\Psi(d_Y(y))=\Psi(mdg_{12}-lx_1x_2x_3)=\Psi(m(a_jf_{ij}-a_if_{ji}-c_{ii}c_{jj})-lx_1x_2x_3)=-lx_1x_2x_3\ne 0$, a contradiction.  Hence $Y$ is not formal.}
 
 \end{example}

  
We give below an example of a finite CW complex  
with only even dimensional cells which is not 
formal. 

\begin{example} {\em 
(i) Let $X=\bs^2\vee \bs^2\vee \bs^2$.  Then $X$ is formal. Computing the  minimal model $\mathcal{M}_X=\Lambda (V)$ 
of $X$ up to degree five we have the following table.
\[ \begin{array}{|l|l|l|l|}
\hline 
 \deg  i &\dim V^i & \textrm{basis} & \textrm{differential} \\ 
\hline
2 & 3 & a_1,a_2, a_3& da_i=0 \\
\hline 
3 & 6 & b_1, b_2, b_3 & db_i=a_i^2\\
&        & b_{12}, b_{23}, b_{13} & db_{ij}=a_ia_j\\
\hline 
4& 6 & c_{ij}, i\ne j& dc_{ij}=b_ia_j-a_ib_{ij}\\
\hline 
5&3  & k_{12}, k_{23}, k_{13}  &dk_{ij}=a_jc_{ij}-a_ic_{ji}+b_ib_j\\
\hline
\end{array}\]

Let $\alpha\in \pi_5(X)$ be such that $\langle k_{12},\alpha\rangle=1$.  Let $Y=X\cup_\alpha e^6$.  Then $Y$ is not formal by Theorem \ref{cell} because the class $[u]\in H^6(Y;\bq)\cong \bq$ is indecomposable. \\
(ii)  Consider the same space $Y$ regarded as a subcomplex of $\wt{Y}=\bc P^3\times \bc P^3\times \bc P^3\cup_\alpha e^6$ where we regard $X$ as the $2$-skeleton of $\wt{X}:=\bc P^3\times \bc P^3\times \bc P^3$.  Note that $\pi_5(\wt{X})=0$ and so $[\alpha]=0$.  It follows that $\wt{Y}$ is formal. Since $Y$ is not formal we see that a subcomplex of a formal space with only even-dimensional cell is not necessarily formal.}
\end{example}

\begin{remark}
(i) By a result of Halperin and Stasheff  \cite[Theorem 1.5]{hs} a nilpotent finite CW complex with only {\it odd} dimensional cells in positive dimensions is formal. Such a CW complex is in fact rationally 
equivalent to a bouquet of odd-dimensional spheres. Halperin and Stasheff point out that this result has also been obtained independently by Baus. \\
(ii) Papadima \cite{papadima} has obtained a criterion for 
the formality of cell attachments.   He also considers spaces $X$ whose cohomology algebra is generated by degree $k$ elements and remarks that formality of such spaces 
can be obtained under the hypothesis that $k\geq c$ (resp. $c-1$) where $c$ is the rational cup-length of $X$ (resp. when $X$ is a Poincar\'e duality space).  
\end{remark}


\subsection{Formality of certain CW complexes}
As an application of Theorem \ref{cell} we shall prove the following theorem.  (This is the corrected version of 
Theorem 1.1 of \cite{cs} which asserted the formality of $X$ without the dimension restriction.)

\begin{theorem} \label{evencw}
Let $X$ be a connected finite CW complex having cells only in dimensions.  If $H^*(X)$ is generated by $H^{2k}(X)$ and $\dim X\le 4k$, then $X$ is formal.
\end{theorem}
\begin{proof}
The $2k$-skeleton $X^{(2k)}$ of $X$ is bouquet of $2k$-dimensional spheres and hence is formal.  We may assume that 
$X$ is of dimension $4k$. 
The minimal 
model $\mathcal{M}=\Lambda V$ of $X^{(2k)}$ has the property that for any standard lower gradation $V=\oplus_{j\ge 0} V_j$, we have $V_0=V^{2k}, V_{\ge 2}\cap V^{4k-1}=0$.  So any element in $\pi_{4-1}(X^{2k})$ is special.
(Note that the Hurewicz homomorphism in dimension $4k-1$ vanishes since $H^{4k-1}(X)=0$.)  It follows that attaching any $4k$-cell to $X^{(2k)}$ 
results in a formal space.   Note that any sub complex of $X$ again has the property that its rational cohomology algebra 
is generated by degree $2k$ elements. We may regard $X=X^{(4k)}$ as obtained from successive cell-attachments 
$X_1,\ldots, X_s$ where $X_{j+1}$ is obtained from $X_{j}$ (with $X_0:=X^{(2k)}$) by  attaching a $4k$-cell, $s$ being the number of $4k$-cells in $X$.  We have just shown that $X_1$ is formal. 

Inductively assume that $X_j$ is formal.  Writing $X_{j+1}=X_j\cup_\alpha e^{4k}$ for a suitable 
$\alpha\in \pi_{4k-1}(X_j)$, we need only show that $\alpha$ is special.  We shall again write $\mathcal{M}=\Lambda V$ 
for the minimal model of $X_j$.  Once again $V_0=V^{2k}$.  We claim that $V_2\cap V^{4k-1}=0$ (and consequently 
$V_j\cap V^{4k-1}=0$ for all $j>2$).   Indeed for dimension reasons $V_1\cap V^p=0$ for all $p<4k-1$.  It follows 
that there are no relations involving elements of $V_0=V^{2k}$  and $V_1\cap V^{4k-1}$ in dimensions less than $6k-1$. 
Thus any element of $\pi_{4k-1}(X_j)$ is special.  It follows by Theorem \ref{cell} that $X_{j+1}$ 
is special.  
\end{proof}

\end{document}